\documentclass{amsart}
\usepackage{amsmath,latexsym,amssymb,amsfonts}

\newtheorem{theorem}{Theorem}[section]
\newtheorem{proposition}[theorem]{Proposition}
\newtheorem{corollary}[theorem]{Corollary}
\newtheorem{lemma}[theorem]{Lemma}

\newtheorem{remark}[theorem]{Remark}

\def\RR{{\mathbb R}}

\begin{document}

\title[The Lebesgue decomposition of the free additive convolution...]{The Lebesgue decomposition of the free additive\\ convolution of two probability distributions}
\author[S. T. Belinschi]{Serban Teodor Belinschi}




\begin{abstract}
We prove that the free additive convolution of two Borel probability 
measures supported on the real line can have a component that is 
singular continuous with respect to the Lebesgue measure on $\mathbb R$
only if one of the two measures is a point mass. The density
of the absolutely continuous part with respect to the Lebesgue measure is shown to 
be analytic wherever positive and finite. 
The atoms of the free additive convolution of Borel probability measures on the real line
have been described by Bercovici and Voiculescu in a previous paper.

\end{abstract}
\maketitle

\vskip 1truecm

\section{Introduction}
The notion of freeness (or free independence) has been introduced by Voiculescu
in \cite{Voiculescu1}, with the main purpose of better understanding free group factors. 
As in the classical case, the distribution of a sum of free random variables is
uniquely determined by the distributions of the summands, and the resulting 
distribution of the sum is called the {\it free additive convolution} of the 
distributions of the summands. 

Specifically, for the case of probability distributions on $\mathbb R$,
let $\mu$ and $\nu$ be Borel probability measures on the real line. We can 
define the free additive convolution $\mu\boxplus\nu$ of $\mu$ and $\nu$ in the 
following way. Denote by $\mathbb F[a,b]$ the free group with
free generators $a$ and $b$, and consider the group von Neumann algebra $L(
\mathbb F[a,b])$ generated by the left regular representation
of $\mathbb F[a,b]$, endowed with the (unique) normal 
faithful trace $\tau$. Choose two selfadjoint operators $X_\mu$ and $X_\nu$ affiliated
with the subalgebras of $L(\mathbb F[a,b])$ generated by the 
images, via the left regular representation, of $a$ and $b$, respectively, so that their 
distribution with respect to $\tau$ is 
$\mu$ and $\nu$, respectively. It has been shown in 
\cite{BVIUMJ} by Bercovici and Voiculescu that the distribution of the selfadjoint 
operator $X_\mu+X_\nu$ with respect to $\tau$ depends only on the distributions 
$\mu$ and $\nu$ of $X_\mu$ and $X_\nu$, respectively. We denote it by $\mu\boxplus\nu$. 
For an introduction to the field of free probability we refer to \cite{VDN}.

An analytic method for the computation of free additive convolutions has been devised 
in \cite{Voiculescu1} (for compactly supported probabilities) and in \cite{BVIUMJ} (for 
the case of probabilities with arbitrary support on $\mathbb R$).  We will give below a 
brief outline of this method.

For any finite positive measure $\sigma$ on $\mathbb R$, define its Cauchy
transform
$$G_\sigma(z)=\int_\mathbb R \frac{d\sigma(t)}{z-t},\quad z\in\mathbb C
\setminus\mathbb R,$$
and let $F_\sigma(z)=1/G_\sigma(z).$
Since $G_\sigma(\overline{z})=\overline{G_\sigma(z)},$ we shall consider from
now on only the restrictions of $F_\sigma$ and $G_\sigma$ to the upper half-plane
$\mathbb C^+=\{z\in\mathbb C\colon \Im z>0\}.$

For given $\alpha\geq0,\beta>0,$ let us denote
$\Gamma_{\alpha,\beta}=\{z\in\mathbb C^+\colon\Im z>\alpha,|\Re z|<\beta\Im z
\}.$ The following two results appear in \cite{BVIUMJ}.
\begin{proposition}\label{inversion+}
Let $\mu$ be a probability on $\mathbb R$. There exists a nonempty domain
$\Omega$ in $\mathbb C^+$ of the form $\Omega=\cup_{\alpha>0}\Gamma_{\alpha,
\beta_\alpha}$ such that $F_\mu$ has a right inverse with respect to
composition $F_\mu^{-1}$ defined on $\Omega$. In addition, we have
$\Im F_\mu^{-1}(z)\leq\Im z$ and $$\lim_{z\to\infty,z\in
\Gamma_{\alpha,\beta}}\frac{F_\mu^{-1}(z)}{z}=1$$ for every $\alpha,\beta>0.$
\end{proposition}
Let $\phi_\mu(z)=F_\mu^{-1}(z)-z$, $z\in\Omega.$ The basic property of the
function $\phi_\mu$ is described in the following theorem of Voiculescu:
\begin{theorem}\label{R-transform}
Let $\mu,\nu$ be two probability measures supported on the real line.
Then $\phi_{\mu\boxplus\nu}(z)=\phi_\mu(z)+\phi_\nu(z)$ for $z$ in the
common domain of the three functions.
\end{theorem}
Thus, the map $\phi$, called the Voiculescu transform, is the free analogue of the logarithm of the 
Fourier transform from classical probability theory.  This function is related to the 
$R$-transform by the equality $\phi_\mu(z)=R_\mu(1/z).$ (Historically, the $R
$-transform, defined as $R_\mu(z)=G_\mu^{-1}(z)-(1/z)$ was first introduced by 
Voiculescu in \cite{Voiculescu1}, but the analysis in the context of measures with 
unbounded support turns out to be simpler when expressed in terms of the Voiculescu transform.)

Another important property for (Cauchy transforms of) free convolutions of probability measures is
subordination. 
It has been shown that $G_{\mu\boxplus\nu}$ is subordinated to $G_\mu$, in the sense that
there exists a unique analytic self-map $\omega$ of the upper half-plane $\mathbb C^+$ so
that $G_{\mu\boxplus\nu}(z)=G_\mu(\omega(z))$, $z\in\mathbb C^+,$ and $\lim_{y\to+\infty}
\omega(iy)/iy=1.$ This result was proved first in \cite{V3} under a genericity 
assumption, then extended to full generality in \cite{Biane1}. A new proof based on 
the theory of fixed points of analytic self-maps of the upper half-plane has been 
given in \cite{Subord}. 

Subordination has been until now the most powerful tool for proving regularity 
results for free convolutions. Pioneering work in this direction has been done by 
Voiculescu, alone in \cite{V3} (see, for example, Proposition 4.7), and together 
with Bercovici in \cite{BVSemigroup} and \cite{BercoviciVoiculescuRegQ}. Among the 
results proved in \cite{BercoviciVoiculescuRegQ},  
we mention the description of the atoms of $\mu\boxplus\nu$ (Theorem 7.4): a number 
$a\in\mathbb R$ is an atom for $\mu\boxplus\nu$ if and only if there exist 
$b,c\in\mathbb R$ so that $a=b+c$ and $\mu(\{b\})+\nu(\{c\})>1.$ Moreover, 
$(\mu\boxplus\nu)(\{a\})=\mu(\{b\})+\nu(\{c\})-1$.
For the special case when $\mu$ is the semicircular distribution (i.e.
$d\mu(t)=\frac{1}{2\pi}\chi_{[-2,2]}(t)\sqrt{4-t^2}dt,$ where $\chi_A$ denotes the 
characteristic function of the set $A$), Biane \cite{Biane2} proved several properties of 
$\mu\boxplus\nu$, from which we mention that the singular continuous part with respect 
to the Lebesgue measure of 
$\mu\boxplus\nu$ is always zero, while the density of its absolutely continuous part
with respect to the Lebesgue measure is bounded and analytic wherever positive. Similar 
results have been proved for measures belonging to partially defined semigroups 
with respect to free additive and multiplicative convolutions (see \cite{BB} and 
\cite{BBercoviciMult}). In \cite{AIHP} it has been shown that, roughly speaking,
 free convolutions of two probability
measures can be purely singular only if at least one of the two measures is a 
point mass. Moreover, for compactly supported measures, the support of the singular part,
if existing, must be of zero Lebesgue measure.

All these results seem to indicate that free convolutions do not favorize large singular 
parts. In Theorem \ref{thm4.1} of this paper we show that when neither $\mu$ nor $\nu$ is a point mass, the 
singular continuous part of $\mu\boxplus\nu$ is zero, while the density of the absolutely
continuous part of $\mu\boxplus\nu$ with respect to the Lebesgue measure is analytic 
outside a closed set of zero Lebesgue measure. 

The rest of the paper is organized as follows: in Section 2 we give without proof several
results from complex analysis that we will use later, in Section 3 we analyze 
the boundary behaviour of the subordination functions, and in Section 4 we prove the
main result of the paper. 

{\bf Acknowledgments.} A weaker result than the one presented in this paper appears 
in the author's PhD thesis. I am deeply grateful to Professor Hari Bercovici for the invaluable help as advisor during 
the  doctoral studies, and afterwards, which made this 
paper possible. My thanks go also to Professor Zhenghan Wang for his help
and support, including through an RAship. Many thanks are due also to the Complex Analysis research 
group 
from Indiana University, especially Professors Eric Bedford, Norm Levenberg, and 
Kevin Pilgrim, for numerous useful discussions.
\section{Preliminary results}

In the following, unless otherwise specified, the attributes "singular", "singular
continuous", and "absolutely continuous" will be considered with respect to the 
Lebesgue measure on $\mathbb R$. Given a finite positive Borel measure $\sigma$ on 
the real line, we denote by $\sigma^s$ (respectively $\sigma^{sc},$ $\sigma^{ac}$)
the singular (respectively singular continuous and absolutely continuous) part of 
$\sigma$. All measures considered in this paper are assumed to be Borel measures.

The following results characterize the Cauchy transform of $\sigma$. For more 
details and proofs we refer to \cite{Achieser}.

\begin{theorem}\label{Cauchy}
Let $G\colon\mathbb C^+\longrightarrow\mathbb C^-$, where $\mathbb C^-=-\mathbb
C^+,$ be an analytic function.
The following statements are equivalent:
\begin{enumerate}
\item There exists a unique positive measure $\sigma$ on $\mathbb R$ such that
$G=G_\sigma$;
\item  For any $\alpha,\beta>0,$ we have that
$$\lim_{z\to\infty,z\in\Gamma_{\alpha,\beta}}zG(z)$$ exists and is finite ($\Gamma_{\alpha,\beta}=\{z\in\mathbb C^+\colon|\Re z|<\alpha\Im z,\Im z>\beta\}$).
\item The limit $\lim_{y\to+\infty}iyG(iy)$ exists and is finite.
\end{enumerate}
Moreover, the limits from 2 and 3 equal $\sigma(\mathbb R).$
\end{theorem}
Observe also that
$$-\frac1\pi\Im G_\sigma(x+iy)=\frac1\pi\int_\mathbb R\frac{y}{(x-t)^2+y^2}d\sigma(t),
\quad x\in\mathbb R,y>0,$$
is the Poisson integral of $\sigma.$

As mentioned also in the introduction, it turns out that in many situations it is much easier to 
deal with the
reciprocal $F_\sigma=1/G_\sigma$ of the Cauchy transform of the measure $\sigma
.$ The following proposition is an obvious consequence of Theorem \ref{Cauchy}:
\begin{proposition}\label{Cauchyalt}
Let $F\colon\mathbb C^+\longrightarrow\mathbb C^+$ be an analytic self-map of
the upper half-plane. The following statements are equivalent:
\begin{enumerate}
\item There exists a positive measure $\sigma$ on $\mathbb R$ such that
$F=1/G_\sigma$;
\item For any $\alpha,\beta>0,$ the limit
$\lim_{z\to\infty,z\in\Gamma_{\alpha,\beta}}\frac{F(z)}{z}$ exists and belongs
to $(0,+\infty)$;
\item The limit $\lim_{y\to+\infty}\frac{F(iy)}{iy}$ exists and belongs to
$(0,+\infty)$.
\end{enumerate}
Moreover, both limits form 2. and 3. equal $\sigma(\mathbb R)^{-1}.$
\end{proposition}
In general, analytic self-maps of the upper half-plane can be represented
uniquely by a triple $(a,b,\rho)$, where $a$ is a real number, $b\in[0,+\infty
),$ and $\rho$ is a positive finite measure on $\mathbb R.$
This representation is called the Nevanlinna representation (see \cite{Achieser}).
\begin{theorem}\label{Nevanlinna}
Let $F\colon\mathbb C^+\longrightarrow\mathbb C^+$ be an analytic function.
Then there exists a triple $(a,b,\rho)$, where $a\in\mathbb R$, $b\geq0,$ and
$\rho$ is a positive finite measure on $\mathbb R$ such that
$$F(z)=a+bz+\int_{\mathbb R}\frac{1+tz}{t-z}d\rho(t),\quad z\in\mathbb C^+.$$
The triple $(a,b,\rho)$ satisfies $a=\Re F(i),$ $b=\lim_{y\to+\infty}\frac{F(iy
)}{iy},$ and $b+\rho(\mathbb R)=\Im F(i)$.
\end{theorem}
The converse of Theorem \ref{Nevanlinna} is obviously true.
\begin{remark}\label{ImF>Imz}
{\rm
An immediate consequence of Proposition \ref{Cauchyalt}  and Theorem \ref{Nevanlinna}
is that for
any finite measure $\sigma$ on $\mathbb R$, we have $\Im F_\sigma(z)\geq
\sigma(\mathbb R)^{-1}\Im z$
for all $z\in\mathbb C^+,$ with equality for any value of $z$ if and only if
$\sigma$ is a point mass. In this case, the measure $\rho$ in the statement of
 Theorem \ref{Nevanlinna} is zero.}
\end{remark}

As observed above, any finite measure $\sigma$ on the real line is uniquely
determined by its
Cauchy transform. Moreover, regularity properties of $\sigma$ can be deduced
from the behaviour of $G_\sigma$, and hence of $F_\sigma$, near the boundary
of its domain. In the following we shall state several classical theorems
concerning analytic self-maps of the unit disc $\mathbb D=\{z\in\mathbb C\colon
|z|<1\}$ and their boundary behaviour, i.e. the behaviour near points belonging
to the boundary $\mathbb T=\{z\in\mathbb C\colon|z|=1\}$ of $\mathbb D$.
Because the upper half-plane is
conformally equivalent to the unit disc via the rational transformation
$z\mapsto\frac{z-i}{z+i},$ most of these theorems will have obvious
formulations for self-maps of the upper half-plane.

We shall consider the set $\mathbb C\cup\{\infty\}$ to be endowed with the usual 
topology: for any point $z\in\mathbb C$, the family $B_n(z)=\{w\in\mathbb C\colon
|z-w|<1/n\}$, $n\in\mathbb N$, forms a basis of neighbourhoods of $z$, while the family
$K_n=\{w\in\mathbb C\colon|w|>n\}\cup\{\infty\}$, $n\in\mathbb N$, forms a basis of 
neighbourhoods for the point infinity. The notions of limit and continuity will be 
considered with respect to this topology and, when subsets of $\mathbb C\cup\{\infty\}
$ are involved, we consider on them the topology inherited from $\mathbb C\cup\{\infty
\}$, unless otherwise specified.
For a function $f\colon\mathbb C^+\longrightarrow\mathbb
C\cup\{\infty\}$,
and a point $x\in\mathbb R$, we say that the nontangential limit of $f$ at
$x$ exists if the limit
$\lim_{z\to x,z\in\Gamma_\alpha(x)}f(z)$ exists in $\mathbb C\cup\{\infty\}$ 
for all $\alpha>0$, where
$\Gamma_\alpha(x)=\{z\in\mathbb C^+\colon|\Re z-x|<\alpha\Im z\}.$
A similar definition holds for functions defined in the unit disc. We shall
denote nontangential limits by $\sphericalangle\lim_{z\to x}f(z),$ or
$$\lim_{\stackrel{ z\longrightarrow x}{{
\sphericalangle}}}f(z).$$
The nontangential limit of $f$ at infinity is defined in a similar way:
$\sphericalangle\lim_{z\to \infty}f(z)$ is said to exist in $\mathbb C\cup\{\infty\}$ 
if the limit $\lim_{z\to\infty,z\in\Gamma_\alpha(0)}f(z)$ exists in $\mathbb C\cup\{\infty\}
$ for all $\alpha>0$.

In the following three theorems are described some properties of meromorphic 
functions in the unit disc related to their nontangential boundary behaviour.

\begin{theorem}\label{Fatou}
Let $f\colon\mathbb D\longrightarrow\mathbb C$ be a bounded analytic function.
Then the set of points $x\in\mathbb T$ at which the nontangential limit of
$f$ fails to exist is of linear measure zero.
\end{theorem}

\begin{theorem}\label{Privalov}
Let $f\colon\mathbb D\longrightarrow\mathbb C$ be an analytic function.
Assume that there exists a set $A$ of nonzero linear measure in $\mathbb T$
such that the nontangential limit of $f$ exists at each point of $A$, and
equals zero. Then $f(z)=0$ for all $z\in\mathbb D$.
\end{theorem}

\begin{theorem}\label{Lindelof}
Let $f\colon\mathbb D\longrightarrow\mathbb C\cup\{\infty\}$ be a meromorphic
function, and let $e^{i\theta}\in\mathbb T$. Assume that the set
$(\mathbb C\cup\{\infty\})\setminus f(\mathbb D)$ contains at least three
points. If there exists a path $\gamma\colon[0,1)\longrightarrow\mathbb D$
such that $\lim_{t\to1}\gamma(t)=e^{i\theta}$ and $\ell=\lim_{t\to1}f(\gamma
(t))$ exists in $\mathbb C\cup\{\infty\}$, then the nontangential limit
of $f$ at $e^{i\theta}$ exists, and equals $\ell$.
\end{theorem}

Theorem \ref{Fatou} is due to Fatou, and Theorem \ref{Privalov} to Privalov.
Theorem \ref{Lindelof} is an extension of a result by Lindel\"{o}f. For
proofs, we refer to \cite{CollingwoodL}, theorems 2.4, 8.1 and 2.20.

These three theorems have obvious reformulations for meromorphic functions
defined on the upper half-plane. For the convenience of the reader we will provide them
below.
\begin{enumerate}
\item Let $f\colon\mathbb C^+\longrightarrow\mathbb C^+\cup\mathbb R\cup\{\infty\}$
be an analytic function. Then the set of points $x\in\mathbb R$ at which the 
nontangential limit of $f$ fails to exist in $\mathbb C\cup\{\infty\}$ is of Lebesgue 
measure zero.
\item Let $f\colon\mathbb C^+\longrightarrow\mathbb C$ be an analytic function. 
Assume that there exists a set $A\subseteq \mathbb R$ of nonzero Lebesgue 
measure such that the nontangential limit of $f$ exists at each point of $A$ and equals
zero. Then $f(z)=0$ for all $z\in\mathbb C^+$.
\item Let $f\colon\mathbb C^+\longrightarrow\mathbb C\cup\{\infty\}$ be a meromorphic
function, and let $x\in\mathbb R\cup\{\infty\}$. Assume that the set
$(\mathbb C\cup\{\infty\})\setminus f(\mathbb C^+)$ contains at least three
points. If there exists a path $\gamma\colon[0,1)\longrightarrow\mathbb C^+$
such that $\lim_{t\to1}\gamma(t)=x$ and $\ell=\lim_{t\to1}f(\gamma
(t))$ exists in $\mathbb C\cup\{\infty\}$, then the nontangential limit
of $f$ at $x$ exists, and equals $\ell$.
\end{enumerate}
The equivalence of statements (2) and (3) above to Theorems \ref{Privalov} and 
\ref{Lindelof} follows by simply considering $f\circ T$, where $T$ is the conformal 
automorphism $T(w)=i\frac{1+w}{1-w},$ and recalling that $T$ preserves angles,
carries a set $A\subseteq\mathbb T$ of zero linear measure into a set 
of zero Lebesgue measure, and vice-versa. The equivalence of Theorem \ref{Fatou}
with statement (1) follows by considering the conjugation of $f$ with $T$ and, 
if necessary, a re-scaling.  We would also like to mention here that in the context of the
above three theorems the point infinity is not in any way a special point; for example,
Theorem \ref{Privalov} forbids any nonconstant analytic function to have constant - 
finite as well as infinite - nontangential limit on a non-negligible set. Thus, given a 
non-constant analytic self-map of $\mathbb C^+$ and a countable set $C\subset
\mathbb C\cup\{\infty\}$, we can always find a set $A\subseteq\mathbb R$ whose 
complement is negligible so that our map has nontangential limits at all 
points of $A$ which do not belong to $C$.

Consider a domain (i.e. an open connected set) $D\subseteq\mathbb C\cup\{\infty
\}$ and a function $f\colon D\longrightarrow\mathbb C\cup\{\infty\}$.
Assume that $\Gamma\subseteq D$ and $x_0$ is an accumulation point for $\Gamma$.
The cluster set $C_\Gamma(f,x_0)$ of the function $f$
at the point $x_0
$ relative to $\Gamma$ is
$$\{z\in\mathbb C\cup\{\infty\}\ |\  \exists
\{z_n\}_{n\in\mathbb N}\subset \Gamma\setminus x_0\ {\rm such\ that }
\lim_{n\to\infty}z_n=x_0, \ \lim_{n\to\infty}f(z_n)=z\}.{\rm }$$
If $\Gamma=D$, we shall write $C(f,x_0)$ instead of $C_D(f,x_0).$
The following result is immediate.
\begin{lemma}\label{ClusterConex}
Let $D\subset\mathbb C\cup\{\infty\}$ be a domain and let
$f\colon D\longrightarrow\mathbb C\cup\{\infty\}$ be continuous. If
$D$ is locally connected at $x\in\overline{D}$, then $C(f, x)$ is connected.
\end{lemma}
This result appears in \cite{CollingwoodL}, as Theorem 1.1.

It will be useful for our purposes to understand the behaviour of analytic self-maps
of $\mathbb C^+$ near open intervals in $\mathbb R$ on which their 
nontangential limits are real almost everywhere.
The following theorem of Seidel can be used to 
describe the behaviour of 
such analytic functions near the boundary of their domain of definition.
For proof, we refer to
\cite{CollingwoodL}, Theorem 5.4.

\begin{theorem}\label{Seidel}
Let $f\colon\mathbb D\longrightarrow\mathbb D$ be an analytic function such
that the radial limit $f(e^{i\theta})=\lim_{r\to1}f(re^{i\theta})$ exists and
has modulus $1$ for almost every
$\theta$ in the interval $(\theta_1,\theta_2)$. If $\theta\in(\theta_1,\theta_2
)$ is such that $f$ does not extend analytically through $e^{i\theta},$ then
$C(f,e^{i\theta})=\overline{\mathbb D}.$
\end{theorem}

This theorem can be applied to self-maps of the upper half-plane, 
via a conformal mapping, but in that case one must consider meromorphic,
instead of analytic, extensions.

A second result refering to the behaviour of $C(f,x)$ for bounded analytic
functions $f$ is the following theorem of Carath\'{e}odory. (This result appears in
\cite{CollingwoodL}, Theorem 5.5.)

\begin{theorem}\label{Caratheodory}
Let $f\colon\mathbb D\longrightarrow\mathbb C$ be a bounded analytic function.
Assume that for almost every $\theta\in(\theta_1,\theta_2)$
the radial limit $f(e^{i\theta})$ belongs to a set $W$ in the plane. Then,
for every $\theta\in(\theta_1,\theta_2)$
the cluster set $C(f, e^{i\theta})$ is contained in the
closed convex hull of $W$.
\end{theorem}

The following proposition is a consequence of the previous two theorems.
\begin{proposition}\label{SeidelCaratheodory}
Let $f$ be an analytic self-map of $\mathbb C^+$ such that
$\lim_{y\to 0}f(x+iy)$ exists and belongs to $\mathbb R$ for almost every
$x\in (a,b).$ Suppose that
$x_0\in(a,b)$ is such that $f$ cannot be continued meromorphically through $x_0
$.
Then for any $c<d$ there is a set
$E\subseteq(a,b)$ of nonzero Lebesgue measure such that
$\lim_{y\to0}f(x+iy)$ exists for all points $x\in E,$ and the set
${\{\lim_{y\to0}f(x+iy)\colon x\in E\}}$ is dense in the interval
$(c,d).$
\end{proposition}
For proof we refer to \cite{AIHP}, Proposition 1.9.

We will next focus on boundary behaviour of derivatives of analytic
self-maps of the unit disk and of the upper half-plane. These results are 
described in detail by Nevanlinna \cite{Nev}  and Shapiro \cite{Shapiro}; see also 
Exercises 6 and 7 in Chapter I of Garnett's book \cite{Garnett}. 

\begin{theorem}\label{JuliaCaratheodory}
Let $f\colon\mathbb D\longrightarrow\mathbb D$ be an analytic function, and 
let $w\in\mathbb T$. The following statements are equivalent:
\begin{enumerate}
\item We have
 $$\liminf_{z\to w}\frac{|f(z)|-1}{|z|-1}<\infty;$$
\item  There exists a number $\zeta\in\overline{\mathbb D}$ such that
$$\lim_{\stackrel{ z\longrightarrow w}{{\sphericalangle}}}f(z)=\zeta,$$
and the limit
\begin{equation}
\ell=
\frac{w}{\zeta}\lim_{\stackrel{ z\longrightarrow w}{{\sphericalangle}}}\frac{f(z)-\zeta}{z-w}
\label{eq1.1}
\end{equation}
exists and belong to $(0,+\infty)$.
\end{enumerate}
Moreover, if the equivalent conditions above are satisfied, the limit
$\sphericalangle\lim_{z\to w}f'(z)$ exists, and the 
following equality holds:
 $$\ell=\frac{w}{\zeta}\lim_{\stackrel{ z\longrightarrow w}{{\sphericalangle}}}
f'(z)=\liminf_{z\to w}\frac{|f(z)|-1}{|z|-1}.$$
If $$\liminf_{z\to w}\frac{|f(z)|-1}{|z|-1}=\infty$$ and 
$$\lim_{\stackrel{ z\longrightarrow w}{{\sphericalangle}}}f(z)=\zeta,$$ then
the limit in equation \eqref{eq1.1} exists and equal infinity.
\end{theorem}
The number $\ell$ from the above theorem is called the Julia-Carath\'{e}odory
derivative of $f$ at $w$.
Since it will be useful in the third and fourth section, we discuss below in some detail the formulation of 
the Julia-Carath\'{e}odory Theorem for self-maps of the upper half-plane. In the following lemma we isolate
the part of the Julia-Carath\'{e}odory Theorem for the upper half-plane that will be used in Sections 3 and 4.
\begin{lemma}\label{JC}
Let $F\colon\mathbb C^+\longrightarrow\mathbb C^+$ be analytic, and let $a\in\mathbb R$. Assume that 
$$\lim_{\stackrel{ z\longrightarrow a}{{\sphericalangle}}}F(z)=c\in\mathbb R.$$ Then
$$\lim_{\stackrel{ z\longrightarrow a}{{\sphericalangle}}}\frac{F(z)-c}{z-a}=
\liminf_{z\to a}\frac{\Im F(z)}{\Im z},$$
where the equality is considered in $\overline{\mathbb C}$. Conversely,
if $$\liminf_{z\to a}\frac{\Im F(z)}{\Im z}<\infty,$$
then $\sphericalangle\lim_{z\to a}F(z)$ exists and belongs to $\mathbb R\cup\{\infty\}$. Moreover, if $F$ is not constant, then we have $\liminf_{z\to a}{\Im F(z)}/{\Im z}>0$.
\end{lemma}
\begin{proof}
Observe that by replacing $F$ with the function $F(z)-c$ we may
assume without loss of generality that $c=0$.
Let $T(z)=\frac{z-i}{z+i},$ $z\in\overline{\mathbb C}$. 
$T$ maps $\mathbb C^+$ conformally onto the unit disc $\mathbb D$, and its 
composition inverse is $T^{-1}(w)=i\frac{1+w}{1-w}.$ 
We consider $f\colon\mathbb D\longrightarrow\mathbb D$ defined by $f(w)=T(F(T^{-1}(w))).$ 
Obviously, $\sphericalangle\lim_{w\to T(a)}f(w)=-1.$ Denote $T(a)=b$.
We have 
\begin{eqnarray*}
\lim_{\stackrel{ z\longrightarrow a}{{\sphericalangle}}}\frac{F(z)}{z-a} & = & \lim_{\stackrel{ w
\longrightarrow b}{{\sphericalangle}}}\frac{F(T^{-1}(w))}{T^{-1}(w)-T^{-1}(b)}\\
& = & \lim_{\stackrel{ w\longrightarrow b}{{\sphericalangle}}}\frac{1+f(w)}{1-f(w)}\cdot\frac{(1-w)(1-b)}{2(w-b)}\\
& = & \frac{(1-b)^2}{4}\lim_{\stackrel{ w\longrightarrow b}{{\sphericalangle}}}\frac{f(w)-(-1)}{w-b}\\
& = & \frac{(1-b)^2}{4}\cdot\frac{-1}{b}\liminf_{w\to b}\frac{1-|f(w)|}{1-|w|}\\
& = & \frac{1}{a^2+1}\liminf_{w\to b}\frac{1-|f(w)|}{1-|w|}.
\end{eqnarray*}
We have used the Julia-Carath\'{e}odory Theorem in the next to last equality and the definition of $b$ in the last 
one. (The situation in which all expressions in the equalities above are infinite is not excluded.)

On the other hand,
\begin{eqnarray*}
\frac{\Im F(z)}{\Im z}=\frac{\Im F(T^{-1}(w))}{\Im T^{-1}(w)} & = & \frac{\Im\left(i\frac{1+f(w)}{1-f(w)}\right)}{\Im\left(
 i\frac{1+w}{1-w}\right)}\\
& = & \frac{\frac{1-|f(w)|^2}{1-2\Re f(w)+|f(w)|^2}}{\frac{1-|w|^2}{1-2\Re w+|w|^2}}\\
& = & \frac{1-|f(w)|}{1-|w|}\cdot\frac{1+|f(w)|}{1+|w|}\cdot
\frac{|1-w|^2}{|1-f(w)|^2}
\end{eqnarray*}
By taking $\liminf$ in the above equality we obtain 
\begin{eqnarray}\label{2}
\liminf_{z\to a}\frac{\Im F(z)}{\Im z} & = & \liminf_{w\to b}\frac{1-|f(w)|}{1-|w|}\cdot\frac{1+|f(w)|}{1+|w|}\cdot
\frac{|1-w|^2}{|1-f(w)|^2}\nonumber\\
& \ge & \liminf_{w\to b}\frac{1-|f(w)|}{1-|w|}\cdot\frac{1+|f(w)|}{1+|w|}\cdot\liminf_{w\to b}
\frac{|1-w|^2}{|1-f(w)|^2}.
\end{eqnarray}
Observe that if $|f(w)|$ does not tend to one, the first $\liminf$ 
in the last row above is infinite. Thus, the first $\liminf$ is
realized on a sequence on which $|f|$ tends to 1. The second $\liminf$ is realized when 
$w$ tends nontangentially to $b$, and equals $
|1-b|^2/4=1/(a^2+1).$ (This follows trivially from the facts that $|f(w)|<1$, $w\in\mathbb 
D$ and $\sphericalangle\lim_{w\to b}f(w)=-1$, since this implies directly that the 
denominator of our expression $\frac{|1-w|^2}{|1-f(w)|^2}$ cannot be greater than $2$,
value reached when $w$ tends to $b$ nontangentially.) Thus,
\begin{eqnarray*}
\liminf_{z\to a}\frac{\Im F(z)}{\Im z} & \ge & \frac{1}{a^2+1}\liminf_{w\to b}\frac{1-|f(w)|}{1-|w|}.
\end{eqnarray*}
We conclude that 
$$\lim_{\stackrel{ z\longrightarrow a}{{\sphericalangle}}}\frac{F(z)}{z-a}\le\liminf_{z\to a}\frac{\Im F(z)}{\Im z}.$$
(Again the case in which both sides of the inequality are infinite is included.) But
$$\lim_{\stackrel{ z\longrightarrow a}{{\sphericalangle}}}\frac{F(z)}{z-a}
=\lim_{y\downarrow0}\frac{F(a+iy)}{iy}
=\lim_{y\downarrow0}\frac{\Re F(a+iy)}{iy}+\frac{\Im F(a+iy)}{y}
\ge\liminf_{z\to a}\frac{\Im F(z)}{\Im z},$$
as our limit is real or infinite.
Thus, 
$$\lim_{\stackrel{ z\longrightarrow a}{{\sphericalangle}}}\frac{F(z)}{z-a}
=\liminf_{z\to a}\frac{\Im F(z)}{\Im z}.$$

Assume now that $\liminf_{z\to a}\Im F(z)/\Im z=d\in\mathbb R_+.$ Equation 
(2) above implies that the limit $\liminf_{w\to b}(1-|f(w)|)/(1-|w|)$ is also finite (recall that $b=T(a)\neq1$, so the 
second $\liminf$ in (\ref{2}) is nonzero),
so, by the Julia-Carath\'{e}odory Theorem $f$ has nontangential limit at
$b$, and thus $F$ has nontangential limit at $a$. (Observe that this limit is infinite 
if and only if $\sphericalangle\lim_{w\to b}f(w)=1$.) Relation (\ref{2}) together with the Julia-Carath\'{e}odory 
Theorem guarantees that $d>0$.
\end{proof}

Consider now an analytic function $f\colon\mathbb D\longrightarrow\overline{
\mathbb D}.$ A point $w\in\overline{\mathbb D}$ is called a Denjoy-Wolff
point for $f$ if one of the following two conditions is satisfied:
\begin{enumerate}
\item[(1)] $|w|<1$ and $f(w)=w$;
\item[(2)] $|w|=1,$ $\sphericalangle\lim_{z\to w}f(z)=w, $ and 
$$\lim_{\stackrel{z\longrightarrow w}{\sphericalangle}}\frac{f(z)-w}{z-w}\leq1.
$$
\end{enumerate}
The following result is due to Denjoy and Wolff.
\begin{theorem}\label{DenjoyWolff}{\large }
Any analytic function $f\colon\mathbb D\longrightarrow\overline{\mathbb D}$ has
a Denjoy-Wolff point. If $f$ has more than one such point, then $f(z)=z$ for all
$z$ in the unit disc. If $z\in\mathbb D$ is a Denjoy-Wolff point for $f$, then
$|f'(z)|\leq1;$ equality occurs only when $f$ is a conformal automorphism of
the unit disc.
\end{theorem}
We refer to Shapiro's book \cite{Shapiro}, Chapter 5, for a detailed introduction to this
subject.

The Denjoy-Wolff point of a function $f$ is characterized also by the fact that
it is the uniform limit on compact subsets of the iterates $f^{\circ n}=
\underbrace{f\circ f\circ\cdots\circ f}_{n\ {\rm times}}$ 
 of $f$. We state the 
following theorem for the sake of completeness (for the original statements,
see \cite{Denjoy} and \cite{Wolff}):

\begin{theorem}
Let $f\colon\mathbb D\longrightarrow\overline{\mathbb D}$ be an analytic
function. If $f$ is not a conformal automorphism of $\mathbb D$, then the 
functions $f^{\circ n}$ converge uniformly on compact subsets of $\mathbb D$ to
the Denjoy-Wolff point of $f$.
\end{theorem}
The previous two theorems have been used in \cite{Subord} to give a new proof for the subordination property for
free additive and free multiplicative convolutions. We reproduce below the result for the free additive convolution
(Theorem 4.1 in \cite{Subord}):
\begin{theorem}\label{sub}
Given two Borel probability measures $\mu,\nu$ on the real line, there exist unique 
analytic functions $\omega_1,
\omega_2\colon\mathbb C^+\longrightarrow\mathbb C^+$ such that 
\begin{enumerate}
\item[{\rm (1)}] $\Im \omega_j(z)\ge\Im z$ for $z\in\mathbb C^+$, and
$$\lim_{y\uparrow+\infty}\frac{\omega_j(iy)}{iy}=1,\quad j=1,2.$$
\item[{\rm (2)}] $F_{\mu\boxplus\nu}(z)=F_\mu(\omega_1(z))=F_\nu(\omega_2(z)),$ and
\item[{\rm (3)}] $\omega_1(z)+\omega_2(z)=z+F_{\mu\boxplus\nu}(z),$ for all $z\in\mathbb C^+.$
\end{enumerate}
\end{theorem}
For $z\in\mathbb C^+,$ the point $\omega_1(z)$ appears as the Denjoy-Wolff point of the function $f_z\colon
\mathbb C^+\longrightarrow\mathbb C^+$ given by $f_z(w)=F_\nu(F_\mu(w)-w+z)-F_\mu(w)+w.$ The function $f_z$ is well-defined for all $z\in\mathbb C^+\cup\mathbb R$
by Remark \ref{ImF>Imz}.
An immediate consequence of Theorem \ref{sub} is the fact that free additive convolution can be defined 
equivalently by purely complex analytic methods, using equations (2) and (3) from the theorem above. This has 
been proved independently by different means in \cite{CG}.

We shall use boundary properties of the subordination functions to describe the atomic,
singular continuous, and 
absolutely continuous parts, with respect to the Lebesgue measure on $\mathbb R
$, of the free convolution $\mu\boxplus\nu$ of two probability measures $\mu,
\nu$ on the real line. 
The following lemma describes the behaviour of the Cauchy transform $G_\mu$
near points belonging to the support of the singular part of the
probability measure $\mu.$ For proof, we refer to \cite{AIHP}, Lemma 1.10, 
\cite{BercoviciVoiculescuRegQ}, Lemma 7.1, and \cite{SteinWeiss}, Theorem 3.16 of
Chapter II.
\begin{lemma}\label{sing}
Let $\mu$ be a Borel probability measure on $\mathbb R$. 
\begin{enumerate}
\item[{\rm(1)}] For $\mu^{s}$-almost 
all $x\in\mathbb R$, the nontangential limit of the Cauchy transform $G_\mu$ of
$\mu$ at $x$ is infinite.
\item[{\rm(2)}] We have $\mu(\{x\})=\sphericalangle\lim_{z\to x}(z-x)G_\mu(z).$
\item[{\rm(3)}] Denote by $f$ the density of $\mu^{ac}$ with respect to the Lebesgue
measure on $\RR$. Then for almost all $x\in\mathbb R$, we have $-\pi f(x)=
\sphericalangle\lim_{z\to x}\Im G_\mu(z).$
\end{enumerate}
\end{lemma}

Finally, we provide a technical lemma, whose proof can be also found in the proof of 
Theorem 2.3 of \cite{AIHP}. We shall give here a more conceptual proof.
\begin{lemma}\label{ntg}
Let $f\colon\mathbb C^+\longrightarrow\mathbb C^+$ be a nonconstant analytic 
function, $x\in
\mathbb R\cup\{\infty\}$, and assume that $C(f,x)\subseteq\mathbb R\cup\{\infty\}$
contains more than one point, and hence, by Lemma \ref{ClusterConex}, a closed
nondegenerate interval or the complement of an open interval. Then for all $c\in C(f,x)$, 
with the possible exception of at most two points, there exists a sequence $\{z_n^{(c)}
\}_{n\in\mathbb N}\subset\mathbb
C^+$ so that 
\begin{enumerate}
\item[{\rm(i)}]$\displaystyle\lim_{n\to\infty}z_n^{(c)}=x$;
\item[{\rm(ii)}]$\displaystyle\lim_{n\to\infty}f(z_n^{(c)})=c$, and
\item[{\rm(iii)}]$\Re f(z_n^{(c)})=c$ for all $n\in\mathbb N$.
\end{enumerate}
\end{lemma}
\begin{proof}
Observe first that, by replacing $f(z)$ with $f\left(-\frac{1}{z}\right)$ if necessary, we 
may assume that $x\in\mathbb R.$ Now pick two arbitrary points $c_1<c_2\in
\mathbb R\cap C(f,x)$, and fix two arbitrary constants $\varepsilon\in(0,|c_1+c_2|/4)$
and $M>\max
\{2+\varepsilon+|c_1|,2+\varepsilon+|c_2|\}.$ Define the compact region 

$$X_M=\{z\in\mathbb C^+\colon|\Re z|\leq M,\ 1/M\leq\Im z\leq M\}.$$
Let $\{z_n^j\}_{n\in\mathbb N}\subset\mathbb C^+$, $j\in\{1,2\}$, be two sequences 
with the following properties:
\begin{enumerate}
\item $\lim_{n\to\infty}z_n^j=x$, $j\in\{1,2\}$;
\item $1>|z_n^2-x|>2|z_n^1-x|>4|z_{n+1}^2-x|$ for all $n\in\mathbb N$;
\item $|f(z_n^j)-c_j|<\frac1n$ for all $n\in\mathbb N$, $n>0$ and $j\in\{1,2\}$.
\end{enumerate}
Define a path $\gamma\colon[0,1]\longrightarrow\mathbb C^+\cup\{x\}$ so that $
\gamma(0)=i,$ $\gamma(1)=x$, $\gamma\left(1-\frac{1}{2n}\right)=z_n^2,$
$\gamma\left(1-\frac{1}{2n+1}\right)=z_n^1,$ and $\gamma$ is linear on the intervals
$\left[0,\frac12\right],$
$\left[1-\frac{1}{2n},1-\frac{1}{2n+1}\right]$ and $\left[1-\frac{1}{2n+1},1-\frac{1}{2n+2}
\right]$ for all $n\in\mathbb N$, $n>0$. We easily observe that $\gamma$ is a simple 
curve in $\mathbb C^+\cup\{x\}$ and $\lim_{t\to1}\gamma(t)=x.$ Thus, there exists $
n_M\in\mathbb N$ so that $f\left(\gamma\left(\left[1-\frac{1}{2n_M},1\right)\right)\right)
\cap X_M=\varnothing.$ Indeed, assume towards contradiction that for all $n\in\mathbb
N$ there exists $t_n\in\left(1-\frac{1}{2n},1\right)$ so that $f(\gamma(t_n))\in X_M
\subset\mathbb C^+$.
Since $X_M$ is compact, there exists a subsequence $\{t_{n_k}\}_k$
of $\{t_n\}_n$ with the property that $\lim_{k\to\infty}f(\gamma(t_{n_k}))$ exists and
belongs to $X_M$. But $t_n\in\left(1-\frac{1}{2n},1\right)$ implies that $\lim_{k\to\infty}
t_{n_k}=1$ and thus $\lim_{k\to\infty}\gamma(t_{n_k})=x$, so that
$\lim_{k\to\infty}f(\gamma(t_{n_k}))\in C(f,x)\subseteq\mathbb R\cup\{\infty\}$. 
Contradiction. Thus indeed the path $f\circ\gamma\colon[0,1)\to\mathbb C^+$
ultimately stays out of $X_M$ as $t\to1$.

For each $n\in\mathbb N$, $n>n_M+\frac1\varepsilon$, define $\Gamma_n$ to be the 
image of the restriction of $f\circ\gamma$ to $\left[1-\frac{1}{2n},1-\frac{1}{2n+1}\right]$.
By the construction of our sequences $\{z_n^j\}_n$ and the path $\gamma$,
$\Gamma_n\subset\mathbb C^+\setminus X_M$ unites through $\mathbb C^+
\setminus X_M$ the balls $B(c_2,1/n)$ and $B(c_1,1/n)$ for all 
$n\in\mathbb N$, $n>n_M+\frac1\varepsilon$. Thus, it must intersect at least one of 
the following three segments:
$$\mathcal S_0=\left\{\frac{c_1+c_2}{2}+is\colon 0<s<1/M\right\},\quad
\mathcal S_1=\{c_1-\varepsilon+is\colon 0<s<1/M\},\quad {\rm or}$$
$$\mathcal S_2=\{c_2+\varepsilon+is\colon 0<s<1/M\},$$
at least once. We conclude that there exists $t_n\in\left(1-\frac{1}{2n},1-\frac{1}{2n+1}
\right)$ so that $f(\gamma(t_n))\in\mathcal S_r$ for some $r=r(n)\in\{0,1,2\}$. Since
this holds for all $n\in\mathbb N$, $n>n_M+\frac1\varepsilon$, there is at least one
of $\mathcal S_0,\mathcal S_1,\mathcal S_2$, call it $\mathcal S_{r_0}$, which is 
intersected by infinitely many paths $\Gamma_n$, and thus there exists a 
subsequence $\{t_{n_k}\}_k$ with $t_{n_k}\in\left(1-\frac{1}{2n_k},1-\frac{1}{2n_k+1}
\right)$ so that $f(\gamma(t_{n_k}))\in\mathcal S_{r_0},$ $k\in\mathbb N$. But
$\lim_{k\to\infty}t_{n_k}=1$, so $\lim_{k\to\infty}\gamma(t_{n_k})=x$, and thus
$\lim_{k\to\infty}f(\gamma(t_{n_k}))=c$, where $\{c\}=\overline{\mathcal S_{r_0}}\cap\mathbb R\subset
\{c_1-\varepsilon,c_2+\varepsilon,(c_1+c_2)/2\}.$

For the point $c$, $\{c\}=\overline{\mathcal S_{r_0}}\cap\mathbb R$, we have constructed
the sequence $z_k^{(c)}=\gamma(t_{n_k})$ satisfying conditions (i), (ii) and (iii) from
our lemma. We shall prove next that all points of $C(f,x)$, with at most two exceptions,
can be realized in the form of a $c$, with $\{c\}=\overline{\mathcal S_{r_0}}\cap\mathbb R$. 
Consider two separate cases:

\noindent{\bf Case 1:} $C(f,x)=\mathbb R\cup\{\infty\}.$ If for all $d\in\mathbb R$
we can find a sequence $\{z_n^{(d)}\}_n$ with the properties (i), (ii) and (iii), then 
we are done. Assume thus that there exists a point  $d\in\mathbb R$ for which no
sequence with
property (iii) exists (observe that properties (i) and (ii) can always be satisfied, by the 
definition of the cluster set). We claim that for any $c\in\mathbb R\setminus\{d\}$ there
is a sequence $\{z_n^{(c)}\}_n$ satisfying properties (i), (ii) and (iii).
Indeed, choose such a point $c\neq d$ and pick $c_1<c<c_2<d$ so that 
$c=(c_1+c_2)/2$. Construct sequences $\{z_n^j\}_n$, $j\in\{1,2\}$,
and a path $\gamma$ as above. We claim that there exists an $\eta>0$ so that
$f(\gamma([0,1)))\cap\{d+is\colon0<s<\eta\}=\varnothing.$ Indeed, assume towards 
contradiction that this is not the case. Then for any $n\in\mathbb N$ there exists
a point $v_n\in f(\gamma([0,1)))\cap\{d+is\colon0<s<1/n\},$ and so there exists
a point $t_n\in[0,1)$ so that $f(\gamma(t_n))=v_n$. Since $\gamma([0,1))\subset
\mathbb C^+$ and $f$ is not constant by hypothesis, the sequence $\{t_n\}_n$
cannot have any accumulation points in $[0,1)$, so $\lim_{n\to\infty}t_n=1$.
But then $\lim_{n\to\infty}\gamma(t_n)=x$, $\lim_{n\to\infty}f(\gamma(t_n))=d$ and
$\Re f(\gamma(t_n))=d$ for all $n\in\mathbb N$. So the sequence $\{\gamma(t_n)\}_n$
satisfies properties (i), (ii) and (iii) for $d$, a contradiction. This proves the existence of 
the required $\eta>0$.

Choose now $\varepsilon\in(0,\min\{\eta,|c_1+c_2|/4,1\})$ and   $M>\max\{1/\eta
,|c_1|+2+\varepsilon,|c_2|+2+\varepsilon\}$ and
define $X_M$ and $\Gamma_n$ as above. It follows immediately then that the 
paths $\Gamma_n$ must intersect the segment
$\mathcal S_0=\left\{\frac{c_1+c_2}{2}+is\colon 0<s<1/M\right\}$ infinitely 
often, since $c_2$ is separated from $c_1$ by $\mathcal S_0$, the lower edge of
$X_M$ and 
$\{d+is\colon0<s<\eta\},$ and $\Gamma_n$ does not intersect the lower edge of
$X_M$ or $\{d+is\colon0<s<\eta\}.$ Thus we have found a sequence $z_n^{(c)}$
as required for our $c=(c_1+c_2)/2$. Since $c\in\mathbb R\setminus\{d\}$ is arbitrary,
we have proved our lemma, with the exceptional points $d$ and $\infty$.

\noindent{\bf Case 2:} $\mathbb R\setminus C(f,x)\neq\varnothing$. This case is
straightforward: pick $c_1$ and $c_2$ to be, one of them the (or a) finite endpoint
of $C(f,x)$, and the other any arbitrary point in  $C(f,x)$ except the previously
chosen one. It follows trivially that any point $c\in C(f,x)$ in between these two points 
accepts a sequence $\{z_n^{(c)}\}_n$ as in our lemma. We leave the details of the 
proof to the reader. Thus, the lemma holds, with the two exceptional points being the 
endpoints of the cluster set.
\end{proof}

\section{Boundary behaviour of the subordination functions}

In the following we fix two Borel probability measures $\mu$ and $\nu$ on $\mathbb R
$, neither of them 
a point mass. For technical reasons, we first study the case when $\mu$ and $\nu$ are 
both convex combinations
of two point masses. (We denote by $\delta_a$ the probability which gives mass one to 
the point $a$.) 
Let $a\in\mathbb R$ be fixed and define the following two self-maps of the upper half-
plane:
$h_\mu(w)=F_\mu(w)-w+a,$ $ h_\nu(w)=F_\nu(w)-w+a,$ $w\in\mathbb C^+.$
The following proposition generalizes Lemma 2.1 of \cite{AIHP}.
\begin{proposition}\label{31}
The function $f_a(w)=h_\nu(h_\mu(w))$ is a conformal automorphism of $\mathbb C^+$ if and only if each of
$\mu$ and $\nu$ are convex combinations of two distinct point masses. In this case, $\omega_1(z)$ and
$\omega_2(z)$, provided by Theorem \ref{sub}, satisfy quadratic equations having $z$ and real numbers
as coefficients.
\end{proposition}
\begin{proof} Assume that $f_a$ is a conformal automorphism of $\mathbb C^+.$ 
We claim that the analytic functions
$h_\mu,h_\nu\colon\mathbb C^+\longrightarrow\mathbb C^+,$ are both
conformal automorphisms of the upper half-plane.
Observe that by the definition of a conformal automorphism, $h_\mu$ must be injective, and $h_\nu$ surjective.
 Let $k$ be the inverse with respect to
composition of $h_\nu\circ h_\mu,$
so that
          $$h_\nu(h_\mu(k(z)))=z,\quad z\in\mathbb C^+.$$
Applying $h_\mu\circ k$ to both sides of the above equality gives
$h_\mu(k(h_\nu(w)))=w$ for all $w$ in the open subset $(h_\mu\circ k)(
\mathbb C^+)$ of the upper half-plane, so, by analytic continuation, for all
$w\in\mathbb C^+.$ This proves surjectivity of $h_\mu$ and injectivity of $h_\nu$ and thus our claim is
proved.

Observe that, by Theorem \ref{Nevanlinna}, $\lim_{y\to+\infty}h_\mu(iy)/iy=\lim_{y\to+\infty}h_\nu(iy)/iy=0.$ 
Thus, there exist real numbers $b_\mu,c_\mu,d_\mu,b_\nu,c_\nu,d_\nu$ so that 
$$h_\mu(z)=\frac{b_\mu z+c_\mu}{z+d_\mu},\quad h_\nu(z)=\frac{b_\nu z+c_\nu}{z+d_\nu},\quad z\in\mathbb
C^+,$$
and
$${\rm det}\left[\begin{array}{cc}
b_\mu & c_\mu \\
1 & d_\mu
\end{array}\right]>0,\quad {\rm det}\left[\begin{array}{cc}
b_\nu & c_\nu \\
1 & d_\nu
\end{array}\right]>0.$$
By the definition of $h_\mu,$ we obtain
\begin{eqnarray*}
F_\mu(z)=h_\mu(z)+z-a & = & \frac{z^2+z(d_\mu+b_\mu-a)+c_\mu-d_\mu a}{z+
d_\mu}\\
& = & \left(t\frac{1}{z-u}+(1-t)\frac{1}{z-v}\right)^{-1},
\end{eqnarray*}
where
$$t=\frac{d_\mu-b_\mu+a+\sqrt{(d_\mu+b_\mu-a)^2-4c_\mu+4d_\mu a}}{2\sqrt{(d_\mu
+b_\mu-a)^2-4c_\mu+4d_\mu a}},$$
and $$u=\frac{a-d_\mu-b_\mu+\sqrt{(d_\mu+b_\mu-a)^2-4c_\mu+4d_\mu a}}{2},$$
$$v=\frac{a-d_\mu-b_\mu-\sqrt{(d_\mu+b_\mu-a)^2-4c_\mu+4d_\mu a}}{2}.$$
Thus, $\mu$ is a convex combination of two point masses $\delta_u$ and $\delta_v,$ with weights $t$ and
$1-t,$ respectively. The result for $\nu$ follows the same way.  
Conversely, if $\mu=t\delta_u+(1-t)\delta_v,$ then a direct computation
shows that $F_\mu(z)=(z-v)(z-u)(z-tv-(1-t)u)^{-1},$ so that
$$h_\mu(z)=\frac{(a-(tu+(1-t)v))z+uv-a(tv+(1-t)u)}{z-(tv+(1-t)u)},\quad z\in\mathbb C^+,$$
and
$${\rm det}\left[\begin{array}{cc}
a-(tu+(1-t)v) & uv-a(tv+(1-t)u) \\
1 & -tv-(1-t)u
\end{array}\right]=t(1-t)(u-v)^2>0,$$
for all $0<t<1,u\neq v$.
This proves the first statement of the
proposition.

Assume now that $\mu=t\delta_u+(1-t)\delta_v,$ $\mu=s\delta_w+(1-s)\delta_x$, with $u\neq v,$ $w\neq x$, and
$0<s,t<1.$ The computations above provide real numbers $b_\mu,c_\mu,d_\mu,b_\nu,c_\nu,d_\nu$ which 
depend polynomially on $u,v,w,x,s,t$ so that 
$$h_\mu(z)=\frac{b_\mu z+c_\mu}{z+d_\mu},\quad h_\nu(z)=\frac{b_\nu z+c_\nu}{z+d_\nu},\quad z\in\mathbb
C^+.$$
On the other hand, by parts (2) and (3) of Theorem \ref{sub} we obtain that 
\begin{eqnarray*}
F_\mu(\omega_1(z)) & = & F_\nu(F_\mu(\omega_1(z))-\omega_1(z)+z)\\
& = & h_\nu(F_\mu(\omega_1(z))-\omega_1(z)+z)+F_\mu(\omega_1(z))-\omega_1(z)+z-a,
\end{eqnarray*}
so that
\begin{eqnarray*}
\omega_1(z) & = & h_\nu(F_\mu(\omega_1(z))-\omega_1(z)+z)+z-a\\
& = & h_\nu(h_\mu(\omega_1(z))-a+z)-a+z\\
& = & \frac{b_\nu\left(\frac{b_\mu\omega_1(z)+c_\mu}{\omega_1(z)+d_\mu}+z-a\right)
+c_\nu}{\frac{b_\mu\omega_1(z)+c_\mu}{\omega_1(z)+d_\mu}+z-a+d_\nu}+z-a\\
& = & \frac{b_\nu(b_\mu \omega_1(z)+c_\mu+(z-a)(\omega_1(z)+d_\mu))+c_\nu
(\omega_1(z)+d_\mu)}{b_\mu\omega_1(z)+c_\mu+(z-a+d_\nu)(\omega_1(z)+d_\mu)}
+z-a.
\end{eqnarray*}
Thus, $\omega_1$ satisfies an equation of degree two with coefficients that are polynomials of $z,s,t,u,v,w,$ and
$x$. The same argument shows the required statement for $\omega_2$
This proves the second statement of the proposition. 
\end{proof}

A consequence of the proposition above is the following
\begin{corollary}\label{2atoms}
With the notations from Theorem \ref{sub}, if $\mu$ and $\nu$ are both convex combinations of two point masses, 
then:
\begin{enumerate}
\item[{\rm(1)}] $\omega_1$ and $\omega_2$ extend continuously to $\mathbb R$ as functions with values in the 
extended complex plane $\overline{\mathbb C}$;
\item[{\rm(2)}] If  either $\omega_1(a)\in\mathbb C^+$ or $\omega_2(a)\in\mathbb C^+$ 
for some $a\in\mathbb R$, then $\omega_1,\omega_2,$ and $F_{\mu
\boxplus\nu}$ extend analytically around $a$;
\item[{\rm(3)}] $(\mu\boxplus\nu)^{sc}=0,$ and $\frac{d(\mu\boxplus\nu)^{ac}(t)}{dt}$ is analytic wherever positive
and finite.
\end{enumerate}
\end{corollary}
Parts (1) and (2) of the corollary are obvious. For the part (3) observe that $F_{\mu
\boxplus\nu}$ extends continuously to $\mathbb R$ and the set 
$F_{\mu\boxplus\nu}^{-1}(\{0\})$ contains only finitely many points, by Theorem
\ref{sub} (3) and by Proposition \ref{31}. Thus, by Lemma \ref{sing} (1), 
$(\mu\boxplus\nu)^{sc}=0$. The statement regarding $\frac{d(\mu\boxplus\nu)^{ac}(t)}{dt}$
follows immediately from Theorem
\ref{sub} (3), Proposition \ref{31}, and Lemma \ref{sing} (3).


 The next theorem describes the boundary behaviour of the subordination functions provided
by Theorem \ref{sub}.

\begin{theorem}\label{Boundary+}
Let $a\in\mathbb R$ be fixed. With the notations from Theorem \ref{sub}, the following hold:
\begin{enumerate}
\item[{\rm (1)}] If either $C(\omega_1,a)\cap\mathbb C^+\neq\varnothing,$ or
$C(\omega_2,a)\cap\mathbb C^+\neq\varnothing,$
then the functions $\omega_1$, $\omega_2,$ and $F_{\mu\boxplus\nu}$ 
extends analytically in a neighbourhood of $a$.
\item[{\rm (2)}] The functions $\omega_1$ and $\omega_2$ have nontangential limits at $a$.
\item[{\rm (3)}] Assume in addition that the sets $I_1=\mathbb R\setminus{\rm supp}
(\mu)$ and $I_2=\mathbb R\setminus{\rm supp}(\nu)$ are nonempty.
Then $\lim_{z\to a}\omega_j(z)$, $j=1,2,$
exist in the extended complex plane $\overline{\mathbb C}.$
\end{enumerate}
\end{theorem}

A slightly less general version of part (3) of this theorem appears implicitly in the proof of 
Theorem 2.3 in \cite{AIHP}.
\begin{proof}
Consider a sequence $\{z_n\}_{n\in\mathbb N}\subset\mathbb C^+$ and a number
$\ell\in\mathbb C^+$ with the property
that $\lim_{n\to\infty}z_n=a$ and $\ell=\lim_{n\to\infty}\omega_1(z_n).$
Define $f\colon(\mathbb C^+\cup\mathbb R)\times\mathbb C^+\longrightarrow\mathbb C^+$ by
$f(z,w)=F_\nu(F_\mu(w)-w+z)-F_\mu(w)+w.$ As noted in the comments following 
Theorem \ref{sub}, $\omega_1(z)$ is the Denjoy-Wolff point of the function $f_z=f(z,\cdot)$ whenever $z\in\mathbb
C^+.$ 
Thus, we
have
$$\ell=\lim_{n\to\infty}\omega_1(z_n)=\lim_{n\to\infty}f(z_n,\omega_1(z_n))=
f(a,\ell).$$
If $f_a=f(a,\cdot)$ is not an automorphism of the upper half-plane, then
we can use Theorem \ref{DenjoyWolff} to conclude that $|f_a'(\ell)|<1.$ It follows from 
Remark \ref{ImF>Imz} 
that $f$ can be extended on a bidisc $B(a,\varepsilon)\times B(\ell,\varepsilon)$ for $
\varepsilon>0$ small 
enough, so by the Implicit Function Theorem there exists $\eta\in(0,\varepsilon)$ and 
an analytic function
$\omega\colon B(a,\eta)\longrightarrow B(\ell,\varepsilon)$ such that $f(z,\omega(z))=
\omega(z)$ for all 
$z\in B(a,\eta)$. Since $f(z,\omega_1(z))=\omega_1(z)$ for $z\in\mathbb C^+$, we 
conclude by the uniqueness 
part of Theorem  \ref{sub} that $\omega$ extends $\omega_1$ to $\mathbb C^+\cup B
(a,\eta).$  Since $F_{\mu\boxplus\nu}=F_\mu\circ\omega_1$ and $\ell=\omega_1(a)\in
\mathbb C^+,$ $F_{\mu\boxplus\nu}$ also extends analytically to some neighbourhood
of $a$. The similar statement for $\omega_2$ follows from part (3) of Theorem 
\ref{sub}. 

The case when $f_a$ is a conformal automorphism of $\mathbb C^+$ is covered by 
Corollary \ref{2atoms}.
This proves (1).


Assume now that the hypothesis of (3) holds and yet $C(\omega_1,a)$
contains more than one point.
By part (1), $C(\omega_1,a)\subseteq
\mathbb R\cup\{\infty\},$ and by Lemma \ref{ClusterConex}, either $C(\omega_1,a)
\setminus\{\infty\}$ is a closed interval in $\mathbb R$ (possibly all of $\mathbb R$), or
$\mathbb R\setminus C(\omega_1,a)$ is an open interval in $\mathbb R.$

By Lemma \ref{ntg},
for any $c$ in $C(\omega_1,a)\setminus\{\infty\},$ with the
possible exception of two points,
there exists a sequence $\{z_n^{(c)}\}_{n\in\mathbb N}$ converging to $a$
such that
$\lim_{n\to\infty}\omega_1(z_n^{(c)})=c,$ and $\Re\omega_1(z_n^{(c)})=c$ for all
$n$ (i.e. $\omega_1(z_n^{(c)})$ converges to $c$ nontangentially - in fact approaches
$c$ vertically).

Thus, if existing, $\lim_{n\to\infty}F_\mu(\omega_1(z_n^{(c)}))=\sphericalangle
\lim_{z\to c}F_\mu(z).$
By Fatou's theorem (Theorem \ref{Fatou}), this limit exists for almost all $c\in C(\omega_1,a).$ Denote it by $F_\mu(c).$
We shall prove that for every $c\in C(\omega_1,a)$ for which $F_\mu(c)$ exists, 
with at most two exceptions, $
F_\mu(c)\not\in\mathbb C^+.$ Indeed,
suppose that $F_\mu(c)\in\mathbb C^+$ for some $c\in C(\omega_1,a)$ for which
$\omega_1(z_n^{(c)})$ as above can be constructed. Recall that
$\omega_1(z_n^{(c)})$ converges nontangentially to $c$. Then,
using 
parts (2) and (3) of Theorem \ref{sub}, Remark
\ref{ImF>Imz}, and the fact that  
$\nu$ is not a point mass,  
we obtain
\begin{eqnarray*}
\Im F_\mu(c) & = & \lim_{n\to\infty}\Im F_\mu(\omega_1(z_n^{(c)}))\\
& = & \lim_{n\to\infty}\Im F_\nu(F_\mu(\omega_1(z_n^{(c)}))-\omega_1(z_{n}^{(c)})+
z_{n}^{(c)})\\
& = & \Im F_\nu(F_\mu(c)-c+a)\\
& > & \Im (F_\mu(c)-c+a)\\
& = & \Im F_\mu(c),
\end{eqnarray*}
which is a contradiction.

 Consider now
${\omega_2}(z)=F_\mu(\omega_1(z))-\omega_1(z)+z,$ $z\in\mathbb C^+.$ We shall
argue that
the set $C(\omega_2,a)\subseteq\mathbb R\cup\{\infty\}$ must be also infinite.
Suppose this were not the case. Then for any $c\in{\rm Int}C(\omega_1,a),$
$$\lim_{z\to a}{\omega_2}(z)=\lim_{n\to\infty}F_\mu(\omega_1(z_n^{(c)}))-
\omega_1(z_n^{(c)})+z_n^{(c)}=F_\mu(c)-c+a,$$
so that, by Theorem \ref{Privalov}, $F_\mu(z)=z-a+\lim_{v\to a}{\omega_2}(v)$
for all $z\in\mathbb C^+$. (We denote by ${\rm Int}A$ the interior of $A\subseteq
\mathbb R$ with respect to the usual topology on $\mathbb R$.)
This contradicts the fact that $\mu$ is not a point mass. So $C(\omega_2,a)$ must be 
an infinite set.

Assume that $c_0\in{\rm Int}C(\omega_1,a)$ is a point where $F_\mu$ does not
continue meromorphically. Proposition \ref{SeidelCaratheodory} shows that the set
$$E=\{c\in C(\omega_1,a)\colon a+F_\mu(c)-c\in I_2\}$$
has nonzero Lebesgue measure.
In particular, for all $c\in E$,
\begin{eqnarray*}
F_\mu(c) & = & \lim_{n\to\infty}F_\mu(\omega_1(z_n^{(c)}))\\
& = & \lim_{n\to\infty}F_\nu(F_\mu(\omega_1(z_n^{(c)}))-\omega_1(z_n^{(c)})+
z_n^{(c)})\\
& = & F_\nu(F_\mu(c)-c+a),
\end{eqnarray*}
where we used the analiticity of $F_\nu$ on $I_2$ in the last equality.
Privalov's theorem (Theorem \ref{Privalov}) implies that $F_\mu(z)=
F_\nu(F_\mu(z)-z+a)$ for all $z\in\mathbb C^+.$
Rewriting this equality gives
                     $$F_\nu(F_\mu(z)-z+a)-(F_\mu(z)-z+a)+a=z,$$
or, equivalently,
$$h_\nu(h_\mu(z))=z,\quad z\in\mathbb C^+.$$
Corollary \ref{2atoms} and Proposition \ref{31} provide now a contradiction.

We conclude that $F_\mu$ extends meromorphically through the whole open interval
${\rm Int}C(\omega_1,a).$  
By the same argument,
$F_\nu$ must continue meromorphically through all of the (nondegenerate) interval
 ${\rm Int}C({\omega_2},a)$.

For any sequence $\{z_n^{(c)}\}_{n\in\mathbb N}$ constructed above, we may take a
subsequence so that
$\lim_{n\to\infty}{\omega_2}(z_n^{(c)})$ exists (we do not claim that ${\omega_2}(
z_n^{(c)})$ converges nontangentially to this limit). Suppose there were a point $d\in C({\omega_2},a)$ and a set $V_{d}
\subset{\rm Int}C(\omega_1,a)$ of nonzero Lebesgue measure
such that $\lim_{n\to\infty}{\omega_2}(z_n^{(c)})=d$ for all $c\in V_{d}.$
Taking limit as $n\to\infty$ in the equality
$$F_\mu(\omega_1(z_n^{(c)}))+z_n^{(c)}=
\omega_1(z_n^{(c)})+{\omega_2}(z_n^{(c)})$$
gives
$$F_\mu(c)+a=c+d,\quad{\rm for\ all}\  c\in V_{d}.$$
Applying again Privalov's theorem, we obtain that $F_\mu(z)=z-(a-d)$ for all
$z\in\mathbb C^+.$ This contradicts the fact that $\mu$ is not a point mass.
Thus, there exists a set $E\subseteq{\rm Int}C(\omega_1,a)$ of positive Lebesgue
measure such that $\{\tilde{c}=\lim_{n\to\infty}{\omega_2}(z_n^{(c)})
\colon c\in E\}\subseteq{\rm Int}C({\omega_2},a).$
Then, since $F_\nu$ extends meromorphically through ${\rm Int}C({\omega_2},a)$,
by Theorem \ref{sub} we conclude that
\begin{eqnarray*}
F_\mu(c) & = &
\lim_{n\to\infty}F_\mu(\omega_1(z_n^{(c)}))\\
& = & \lim_{n\to\infty}F_\nu({\omega_2}(z_n^{(c)}))\\
& = & \lim_{n\to\infty}F_\nu\left({F_\mu(\omega_1
(z_n^{(c)}))}-\omega_1(z_n^{(c)})+z_n^{(c)}\right)\\
& = & F_\nu(a+F_\mu(c)-c)
\end{eqnarray*}
for all $c\in E$. Privalov's theorem implies that $F_\mu(z)=F_\nu(a+F_\mu(z)-
z)$ for all $z\in\mathbb C^+.$ As we have proved already, this implies that
$h_\mu$ and $h_\nu$ are conformal automorphisms of the upper half-plane. Corollary
\ref{2atoms} and Proposition \ref{31} provide again a contradiction.
This proves part (3) of the theorem. 


Consider now the case when at least one of the two sets $I_1,I_2$ is empty.
Without loss of generality, assume $I_1=\varnothing$.
We show that in this case, if any of the two sets $C(\omega_1,a)$, $C(\omega_2,a)$
contains more than one point, we must have $C(\omega_1,a)=C(\omega_2,a)=
\mathbb R\cup\{\infty\}.$ We shall do this in four steps. 
\newline\noindent
{\bf Step 1:} We show that ${\rm Int}C(\omega_1,a)\cap{\rm supp}
(\mu^{ac})=\varnothing.$ Indeed, assume this is not the case. Then (by Lemmas
\ref{sing} and \ref{ntg}) for
almost all $c$ in ${\rm Int}C(\omega_1,a)\cap{\rm supp}
(\mu^{ac})$ with respect to the Lebesgue measure, we have $\sphericalangle
\lim_{z\to c}F_\mu(z)\in\mathbb C^+$, and there exists $\{z_n^{(c)}\}_{n\in
\mathbb N}$ so that $\lim_{n\to\infty}z_n^{(c)}=a$, $\Re\omega_1(z_n^{(c)})
=c,$ and $\lim_{n\to\infty}\omega_1(z_n^{(c)})=c.$ Thus,
\begin{eqnarray*}
\lim_{n\to\infty}\omega_2(z_n^{(c)}) & = & \lim_{n\to\infty}
F_{\mu}(\omega_1(z_n^{(c)}))-\omega_1(z_n^{(c)})+z_n^{(c)}\\
& = & a-c+\lim_{\stackrel{z\longrightarrow c}{{\sphericalangle}}}
F_\mu(z)\in\mathbb C^+,
\end{eqnarray*}
contradicting part (1) of the theorem. 
\newline\noindent
{\bf Step 2:} We show that $C(\omega_2,a)=\mathbb R\cup\{\infty\}.$ One 
inclusion is an immediate consequence of part (1) of this theorem. Thus, it
is enough to show that $C(\omega_2,a)$ is dense in $\mathbb R$. By Step 1 
and our hypothesis $I_1=\varnothing,$ we have ${\rm Int}C(\omega_1,a)
\subseteq{\rm supp}(\mu^{s})\setminus{\rm supp}(\mu^{ac}).$  Thus, by
Proposition \ref{SeidelCaratheodory}, for any $x\in{\rm Int}C(\omega_1,a)$
and any open interval $I$ containing $x$, the set

$$\{\lim_{\stackrel{z\longrightarrow c}{{\sphericalangle}}}
F_\mu(z)\colon c\in I, F_\mu\ {\rm has\ nontangential\ limit\ at}\ c\}$$
is dense in $\mathbb R$. So, given $x$ as above, $\varepsilon>0$ and 
$s\in\mathbb R$, there exists $c\in I\cap(x-\frac{\varepsilon}{2},x+
\frac{\varepsilon}{2})$ such that $|s-x+a-\sphericalangle\lim_{z\to c}F_\mu(z)|
<\varepsilon/2,$ and thus
\begin{eqnarray*}
|s-\lim_{n\to\infty}\omega_2(z_n^{(c)})| & = & |s-\lim_{n\to\infty}F_{\mu}
(\omega_1(z_n^{(c)}))-\omega_1(z_n^{(c)})+z_n^{(c)}|\\
& = & |s-c+a-\lim_{\stackrel{z\longrightarrow c}{{\sphericalangle}}}
F_\mu(z)|\\
& \le & |s-x+a-\lim_{\stackrel{z\longrightarrow c}{{\sphericalangle}}}
F_\mu(z)|+|x-c|\\
& < & \frac{\varepsilon}{2}+\frac{\varepsilon}{2}\\
& = & \varepsilon.
\end{eqnarray*}
We conclude that $C(\omega_2,a)$ is dense in $\mathbb R$. This proves step 2.
\newline\noindent
{\bf Step 3:} We show that $\nu=\nu^s.$ This follows from Step 2 and the argument 
used in Step 1.
\newline\noindent
{\bf Step 4:} We show that $C(\omega_1,a)=\mathbb R\cup\{\infty\}.$ If $\nu$ is 
a convex combination of point masses and there exist two consecutive atoms of 
$\nu$ at $\alpha$ and $\beta$ ($\alpha<\beta$), then obviously $G_\nu$ extends 
analytically to the interval $(\alpha,\beta)$ and $G_\nu((\alpha,\beta))=\mathbb R$,
so that $F_\nu((\alpha,\beta))=(\mathbb R\setminus\{0\})\cup\{\infty\}.$ 
For any $x\in(\alpha,\beta)$ let $z_n^{(x)}$ converge to $a$ so that 
$\omega_2(z_n^{(x)})\to x$ as $n$ tends to infinity. Then
\begin{eqnarray*}
\lim_{n\to\infty}\omega_1(z_n^{(x)}) & = & \lim_{n\to\infty}F_\nu(\omega_2(z_n^{(x)}))
-\omega_2(z_n^{(x)})+z_n^{(x)}\\
& = & F_\nu(x)-x+a,
\end{eqnarray*}
and thus
$$C(\omega_1,a)\supseteq\overline{\{\lim_{n\to\infty}\omega_1(z_n^{(x)})\colon 
x\in(\alpha,\beta)\}}=\overline{\{F_\nu(x)-x+a\colon x\in(\alpha,\beta)\}}=\mathbb R
\cup\{\infty\}.$$

If either $\nu$ is not purely atomic, or there exist no consecutive atoms of $\nu$, 
then there exists at least one point $x_0\in\mathbb R$ so that $F_\nu$ does not 
extend meromorphically through $x_0$. The argument used in the proof of Step 2, 
with $I$ an open interval containing $x_0$, assures us that $C(\omega_1,a)=
\mathbb R\cup\{\infty\}.$
This proves Step 4.

We have proved now that if there exists a point $a$ where either $C(\omega_1,a)$ or
$C(\omega_2,a)$ is nondegenerate (i.e. contains more than one point), then
$\mu=\mu^s$, $\nu=\nu^s$, at least one has total support, and $C(\omega_1,a)=
C(\omega_2,a)=\mathbb R\cup\{\infty\}.$ We assume without loss of generality that $\mu$ has support equal to the 
whole real line. Choose a point $c\in\mathbb R$ so that 
$\sphericalangle\lim_{z\to c}F_\mu(z)=0$ and there exists $\{z_n^{(c)}
\}_{n\in\mathbb N}$ converging to $a$ so that $\Re\omega_1(z_n^{(c)})=c$ and
$\lim_{n\to\infty}\omega_1(z_n^{(c)})=c.$ By Lemmas \ref{sing} and \ref{JC} we have 
\begin{eqnarray*}
\frac{1}{\mu(\{c\})}-1 & = & \lim_{\stackrel{z\longrightarrow c}{{\sphericalangle}}}
\frac{F_\mu(z)}{z-c}-1\\
& = & \liminf_{z\to c}\frac{\Im F_\mu(z)}{\Im z}-1\\
& \le & \liminf_{n\to\infty}\frac{\Im F_\mu(\omega_1(z_n^{(c)}))}{\Im 
\omega_1(z_n^{(c)})}+\frac{\Im z_n^{(c)}}{\Im\omega_1(z_n^{(c)})}-1\\
& = & \liminf_{n\to\infty}\frac{\Im \omega_2(z_n^{(c)})}{\Im\omega_1(z_n^{(c)})}\\
& = & \left(\limsup_{n\to\infty}
\frac{\Im \omega_1(z_n^{(c)})}{\Im\omega_2(z_n^{(c)})}\right)^{-1}\\
& = & \left(\limsup_{n\to\infty}\frac{\Im F_\nu(\omega_2(z_n^{(c)}))}{\Im 
\omega_2(z_n^{(c)})}+\frac{\Im z_n^{(c)}}{\Im\omega_2(z_n^{(c)})}-1\right)^{-1}\\
& \le & \left(\liminf_{n\to\infty}\frac{\Im F_\nu(\omega_2(z_n^{(c)}))}{\Im 
\omega_2(z_n^{(c)})}-1\right)^{-1}\\
& \le & \left(\liminf_{z\to a-c}\frac{\Im F_\nu(z)}{\Im z}-1\right)^{-1}.
\end{eqnarray*}
By our assumption, $(\mu(\{c\}))^{-1}-1\in(0,\infty]$, so that 
$$\liminf_{z\to a-c}\frac{\Im F_\nu(z)}{\Im z}<\infty.$$ Lemma 
\ref{JC} 
implies that $F_\nu$ has a nontangential limit at $a-c$ belonging to $\mathbb R\cup\{\infty\}$. Denote it by
$l.$ Moreover, we claim 
that $\mu(\{c\})>0.$ Indeed, assume this is not the case. Then the above chain
of inequalities together with Lemma \ref{JC}   
imply that 
$$\liminf_{z\to a-c}\frac{\Im F_\nu(z)}{\Im z}=1,$$
so that $h_\nu(z)=F_\nu(z)-z$ is constant, by Lemma \ref{JC}.
This contradicts the assumption that $\nu$ is not a point mass and proves
our claim. We conclude that $\mu$ must be an infinite convex combination of point
masses, densely distributed in $\mathbb R$.

Assume towards contradiction that $\omega_1$ has no nontangential limit at $a$. Then
there exists an angle $\Gamma\subset\mathbb C^+$ with vertex at $a$ and bisected
by the line $a+i\mathbb R_+$ with the property that $C_\Gamma(\omega_1,a)$ is
infinite.
By Lemma \ref{ClusterConex} and the argument above, $C_\Gamma(\omega_1,a)$ 
must 
contain a nondegenerate open subinterval $J$ so that for any $c\in J$ there exists
a sequence $\{z_n^{(c)}\}_{n\in\mathbb N}\subset\Gamma$ converging to $a$ such
that $\Re \omega_1(z_n^{(c)})=c$ and $\Im\omega_1(z_n^{(c)})$ converges to zero 
as $n\to\infty$. Let $m=\max_{x\in\mathbb R}\nu(\{x\}).$ By hypothesis, $0\le m<1.$
Since the set of atoms of $\mu$ is dense in $\mathbb R$, in particular infinite, there
 exists $c\in J$ 
so that $0<\mu(\{c\})<1-m$. 
Observe that since $1>\mu(\{c\})>0,$ the non-constant function $F_\mu(z)-z$ maps any nontangential path ending
at $c$ into a nontangential path. Indeed, since
$$
\lim_{\stackrel{z\longrightarrow c}{{\sphericalangle}}}\frac{(F_\mu(z)-z)+c}{z-c}=
\frac{1}{\mu(\{c\})}-1\in(0,+\infty),$$
this follows by simply analyzing the real and the imaginary part of the left-hand term 
above.
Thus, the sequence $\{
\omega_2(z_n^{(c)})\}_{n\in\mathbb N}=\{F_{\mu}(\omega_1(z_n^{(c)}))
-\omega_1(z_n^{(c)})+z_n^{(c)}\}_{n\in\mathbb N}$  must also converge
nontangentially to $a-c.$ But we have seen above that $F_\nu$ has nontangential
limit at $a-c$. Thus,
$$0=\lim_{\stackrel{z\longrightarrow c}{{\sphericalangle}}}F_\mu(z)=
\lim_{n\to\infty}F_\mu(\omega_1(z_n^{(c)}))=\lim_{n\to\infty}F_\nu(
\omega_2(z_n^{(c)}))=\lim_{\stackrel{z\longrightarrow a-c}{{\sphericalangle}}}F_\nu(z),$$
and
\begin{eqnarray*}
\frac{1}{\mu(\{c\})}-1 & \leq & \left(\liminf_{z\to a-c}
\frac{\Im F_\nu(z)}{\Im z}-1\right)^{-1}\\
& = & \left(\lim_{\stackrel{z\longrightarrow a-c}{{\sphericalangle}}}
\frac{F_\nu(z)}{z-(a-c)}-1\right)^{-1}\\
& = & \left(\frac{1}{\nu(\{a-c\})}-1\right)^{-1}.
\end{eqnarray*}
Multiplication with $\frac{1}{\nu(\{a-c\})}-1$ in the above inequality gives
$$(1-\mu(\{c\}))(1-\nu(\{a-c\}))\leq\mu(\{c\})\nu(\{a-c\}),$$ so
$$\mu(\{c\})+\nu(\{a-c\})\ge1.$$ This contradicts the choice $\mu(\{c\})<1-m$ 
and concludes the proof of part (2) of the theorem.
\end{proof}


\section{The main result}

We can now prove the main result of this paper. For the sake of completeness, we
state in the theorem below the result of Bercovici and Voiculescu describing the
atoms of the free additive convolution of two probability distributions.
\begin{theorem}\label{thm4.1}
Let $\mu,\nu$ be two Borel probability measures on $\mathbb R$, neither of them
a point mass. Then
\begin{enumerate}
\item[{\rm(1)}] The point $a\in\mathbb R$ is an atom of the measure $\mu
\boxplus\nu$ if and only if there exist $b,c\in\mathbb R$ such that
$a=b+c$ and $\mu(\{b\})+\nu(\{c\})>1.$ Moreover, $(\mu\boxplus\nu)(\{a\})=
\mu(\{b\})+\nu(\{c\})-1.$ 
\item[{\rm (2)}] The absolutely continuous part of $\mu\boxplus\nu$ 
is always nonzero, and its density is analytic wherever positive and finite. More 
precisely, there exists an open set $U\subseteq\mathbb R$ so that the density function 
$f(x)=\frac{d(\mu\boxplus\nu)^{ac}(x)}{dx}$ with respect to the Lebesgue measure in the real line is analytic on $U$ and $(\mu\boxplus\nu)^{ac}(\mathbb R)=\int_{U}f(x)dx$.
\item[{\rm(3)}] The singular continuous part of $\mu\boxplus\nu$ is zero.
\end{enumerate}
\end{theorem}
\begin{proof}
Part (1) of the theorem is due to Bercovici and Voiculescu (see 
\cite{BercoviciVoiculescuRegQ}, Theorem 7.4).
We shall proceed with the proof of part (2).
Suppose that $\mu\boxplus\nu$ is purely singular, and thus for almost
all $x\in\mathbb R$ with respect to the Lebesgue measure, we have
$$\lim_{y\to0}\Im F_{\mu\boxplus\nu}(x+iy)=\lim_{y\to0}\Im 
G_{\mu\boxplus\nu}(x+iy)=0.$$ 
By part (1), we are assured that $\mu\boxplus\nu$ cannot be purely atomic, so 
we must have $(\mu\boxplus\nu)^{sc}\neq0,$
and hence, by Lemma \ref{sing}, $\lim_{y\to0} F_{\mu\boxplus\nu}(x+iy)=0$ 
for uncountably many $x\in\mathbb R.$
Theorem \ref{Seidel} applied to the function $F_{\mu\boxplus\nu}$ yields
a point $x_0\in\mathbb R$ such that $C(F_{\mu\boxplus\nu},x_0)=
\overline{\mathbb C^+}.$ Using relations (2) and (3) from Theorem \ref{sub} we conclude that at 
least one of 
$C(\omega_1,x_0),C(\omega_2,x_0)$, will intersect the upper
half-plane.
But now we can apply Theorem \ref{Boundary+} (1)  
to 
obtain a contradiction. 
Thus, $\mu\boxplus\nu$ cannot be purely singular\footnote{This is a 
slightly modified version of part of the proof of Corollary 2.5 in
\cite{AIHP}. Of course, the result presented here is much stronger.}.

Next we prove that there exists an open subset $U$ of $\mathbb R$ on which
the density $f(x)=\frac{d(\mu\boxplus\nu)^{ac}(x)}{dx}$ is analytic. 
By Theorem \ref{Fatou}, there exists a subset $E$ of $\mathbb R$ of zero 
Lebesgue measure such that for all $x\in\mathbb R\setminus E$ the limits
$\lim_{y\to0}F_{\mu\boxplus\nu}(x+iy)$, $\lim_{y\to0}\omega_j(x+iy),$ 
$j\in\{1,2\}$ exist and are finite. Also, by Lemma \ref{sing}, (3), for almost all 
$x\in{\rm supp}((\mu\boxplus\nu)^{ac})\setminus E$, with respect to
 $(\mu\boxplus\nu)^{ac},$ we have $\lim_{y\to0}F_{\mu\boxplus\nu}(x+iy)\in
 \mathbb C^+.$ By relation (3) in Theorem \ref{sub}, at 
least one of $\lim_{y\to0}\omega_j(x+iy),$ $j\in\{1,2\},$ 
must also be in $\mathbb C^+$. For definiteness,
assume $\omega_1(x)=\lim_{y\to0}\omega_1(x+iy)\in\mathbb C^+.$ Part (1) of
Theorem \ref{Boundary+}  assures us that
$\omega_1,\omega_2$ and $F_{\mu\boxplus\nu}$ extends analytically through
$x$. We conclude by part (3) of Lemma \ref{sing} that the density of 
$(\mu\boxplus\nu)^{ac}$ must be analytic in $x$. 
 Of course,
the set $U$ of points $x$
where $f(x)=\frac{d(\mu\boxplus\nu)^{ac}(x)}{dx}$ is analytic is open in 
$\mathbb R$, and by Lemma \ref{sing} (3) 
we conclude that
$\int_Uf(x)dx=(\mu\boxplus\nu)^{ac}(\mathbb R)$.

To prove (3), assume that $(\mu\boxplus\nu)^{sc}\neq0;$ Lemma \ref{sing}
provides an uncountable set $H$ of points $c\in\mathbb R$ such that 
$(\mu\boxplus\nu)(\{c\})=0$ and $\sphericalangle\lim_{z\to c}
F_{\mu\boxplus\nu}(z)=0.$
Theorem \ref{Boundary+} (2) assures us that $\omega_1$ and $\omega_2$ have 
nontangential limits at all points of $\mathbb R$, in particular at
each such point $c\in H$. We claim that
\begin{enumerate}
\item[{(i)}] $\sphericalangle\lim_{z\to c}\omega_j(z)\in\mathbb R$, $j=1,2$. Denote 
those limits by $v_j$, $j=1,2$;
\item[{(ii)}] The following equalities hold:
$$\lim_{\stackrel{ z\longrightarrow v_1}{{\sphericalangle}}}F_\mu(z)=0
\quad {\rm and}\quad
\lim_{\stackrel{ z\longrightarrow v_2}{{\sphericalangle}}}F_\nu(z)=0.$$
\item[{(iii)}] $\mu(\{v_1\})+\nu(\{v_2\})=1.$
\end{enumerate}
Parts (2) and (1) of Theorem \ref{Boundary+} guarantee the existence of the
limits in (i) and the fact that they do not belong to $\mathbb C^+.$
If, say, $\sphericalangle\lim_{z\to c}\omega_1(z)=\infty$, then Theorem \ref{Lindelof} 
assures us that
\begin{eqnarray*}
0 & = & \lim_{y\downarrow0}F_{\mu\boxplus\nu}(c+iy)\\
& = & \lim_{y\downarrow0}F_\mu
(\omega_1(c+iy))\\
& = & \lim_{s\to\infty,s\in\omega_1(c+i\mathbb R)}F_{\mu}(s)\\
& = & \lim_{\stackrel{ z\longrightarrow \infty}{{\sphericalangle}}}F_\mu(z)\\
& = & \infty,
\end{eqnarray*}
which is of course a contradiction.
Thus (i) holds.
To prove (ii) just observe that, by the same Theorem \ref{Lindelof}
\begin{eqnarray*}
0 & = & \lim_{\stackrel{ z\longrightarrow c}{{\sphericalangle}}}F_{\mu\boxplus\nu}(z)\\
& = &\lim_{\stackrel{ z\longrightarrow c}{{\sphericalangle}}}F_\nu(\omega_2(z))\\
& = & \lim_{z\to v_2,z\in\omega_2(c+i\mathbb R)}F_\nu(z)\\
& = & \lim_{\stackrel{ z\longrightarrow v_2}{{\sphericalangle}}}F_\nu(z).
\end{eqnarray*}
The same argumend provides the proof for $F_\mu$.

We now prove (iii). Recall that by Lemmas \ref{sing} and \ref{JC}, for any $c$ as above
we have 
\begin{eqnarray*}
\frac{1}{\mu(\{c\})} & = & \lim_{\stackrel{z\longrightarrow c}{\sphericalangle}}\left[(z-c)
\int_\mathbb R\frac{d\mu(t)}{z-t}\right]^{-1}\\
& = & \lim_{\stackrel{z\longrightarrow c}{\sphericalangle}}\frac{F_\mu(z)}{z-c}\\
& = & \liminf_{z\to c}\frac{\Im F_\mu(z)}{\Im z}.
\end{eqnarray*}
Thus the following chain of inequalities holds:
\begin{eqnarray*}
\frac{1}{\nu(\{v_2\})}-1 & = & \liminf_{z\to v_2}\frac{\Im F_\nu(z)}{\Im z}-1\\
& \leq & \liminf_{y\downarrow0}\frac{\Im
F_\nu(\omega_2(c+iy))}{\Im\omega_2(c+iy)}-1\\
& \leq & \limsup_{y\downarrow0}\left(\frac{\Im F_\nu(\omega_2(c+iy))}{\Im
\omega_2(c+iy)}+
\frac{y}{\Im \omega_2(c+iy)}-1\right)\\
& = & \limsup_{y\downarrow0}\frac{\Im \omega_1(c+iy)}{\Im \omega_2(c+iy)}\\
& = & \left(\liminf_{y\downarrow0}\frac{\Im\omega_2(c+iy)}{\Im \omega_1(c+iy)}
\right)^{-1}\\
& = & \left(\liminf_{y\downarrow0}\left(\frac{\Im F_\mu(\omega_1(c+iy))}{\Im
\omega_1(c+iy)}+\frac{y}{\Im \omega_1(c+iy)}\right)-1\right)^{-1}\\
& \leq & \left(\liminf_{y\downarrow0}\frac{\Im F_\mu(\omega_1(c+iy))}{\Im \omega_1(c+
iy)}+
\liminf_{y\downarrow0}\frac{y}{\Im \omega_1(c+iy)}-1\right)^{-1}\\
& \leq & \left(\liminf_{y\downarrow0}\frac{\Im F_\mu(\omega_1(c+iy))}{\Im \omega_1(c
+iy)}-1\right)^{-1}\\
& \leq & \left(\liminf_{z\to v_1}\frac{\Im F_\mu(z)}{\Im z}-1\right)^{-1}\\ 
& = & \left(\frac{1}{\mu(\{v_1\})}-1\right)^{-1}.
\end{eqnarray*}
We have assumed that $\mu$ and $\nu$ are not point masses, so the above 
implies that $1<\frac{1}{\mu(\{v_1\})},\frac{1}{\nu(\{v_2\})}<
\infty.$ Thus, multiplying the above inequality by 
$\frac{1}{\mu(\{v_1\})}-1$ will give 
$$(1-\mu(\{v_1\}))(1-\nu(\{v_2\}))\leq\mu(\{v_1\})\nu(
\{v_2\}),$$ or, equivalently,
$$\mu(\{v_1\})+\nu(\{v_2\})\geq1.$$
Using relation (3) from Theorem \ref{sub} and the fact that $F_{\mu\boxplus\nu}(c)=0$,
we obtain $$v_1+v_2=c\quad {\rm for\ all}\ c\in H.$$
Since $c$ has been chosen so that $(\mu\boxplus\nu)(\{c\})=0$, part (1) of the theorem 
tells us that the inequality above must be an equality. This proves the last point of our
claim.

Now, since any probability can have at most countably many atoms, this, together 
with part ({iii}) of our claim contradicts the fact that $H$ is uncountable
and concludes the proof.
\end{proof}

MSC: primary, 46L54; secondary, 30D40.

Address:

\noindent Department of Pure Mathematics, University of Waterloo, 
              200 University Street West, Waterloo, ON, N2L 1V9,
              Canada,\\
              and\\
              Institute of Mathematics "Simion Stoilow" of the Romanian Academy,
              Bucharest, Romania.\\
              Tel.: +1-(519)-888-4567 Ext. 32803\\
              E-mail: {sbelinsc@math.uwaterloo.ca}            \\
             \emph{Present address:Department of Pure Mathematics, University of Waterloo, 
              200 University Street West, Waterloo, ON, N2L 1V9,
              Canada\\}


\begin{thebibliography}{3}

\bibitem{Achieser} {Akhieser, N. I.}: { 
The classical moment problem and some related questions in analysis}.
{Hafner Publishing Co.}, {New York}, {1965}




\bibitem{AIHP}{Belinschi, S. T.}: { A note on regularity for free convolutions}. Ann. Inst. H. Poincar\'{e} Probab. Statist. {\bf 42}, {no. 3}, {635--648} {(2006)}.

\bibitem{BB} {Belinschi, S. T., Bercovici, H}: {Atoms and regularity for
measures in a
partially defined free convolution semigroup}. Math. Z., {\bf 248}, {no. 4}, {665--674}
{(2004)}.

\bibitem{BBercoviciMult} {Belinschi, S. T., Bercovici, H.}: {
Partially Defined Semigroups Relative to Multiplicative Free Convolution}.
{Internat. Math. Res. Not.}, {no.2},
{65--101} {(2005)}.


\bibitem{Subord}{Belinschi, S. T., Bercovici, H.}: {A new approach to subordination results
in free probability}. {J. Anal. Math.}, {to appear}.



\bibitem{BVIUMJ} {Bercovici, H., Voiculescu, D.}: {
Free convolutions of measures with unbounded support}.
{Indiana Univ. Math. J.}, {\bf 42}, no.3, {733--773} {(1993)}.

\bibitem{BVSemigroup} {Bercovici, H., Voiculescu, D}: {Superconvergence to the central limit
and failure of the Cram$\acute{e}$r theorem for free random variables}.
{Probab. Theory Related Fields} {\bf 102},  {215--222} {(1995)}.


\bibitem{BercoviciVoiculescuRegQ} {Bercovici, H., Voiculescu, D.}: {Regularity questions for free convolution},
{Nonselfadjoint operator algebras, operator theory, and related topics, 37--47,
Oper. Theory Adv. Appl.} {104}, {Birkh\"{a}user, Basel}, {1998}.

\bibitem{Biane1} {Biane, P.}: {Processes with free increments}. {Math. Z.}
{\bf 227}, {143--174} {(1998)}.

\bibitem{Biane2}{Biane, P}: {On the Free Convolution with a Semi-circular Distribution}.
Indiana Univ. Math. J. {\bf46}, 705--718 (1997).


\bibitem{C} {Carath\'eodory, C}: {\"Uber die Winkelderivierten von
beschr\"ankten analytischen
Funktionen}. Sitzungbericht Preuss. Akad. Phys. Math. {\bf 4}, 1--18 (1929).

\bibitem{CG} Chistyakov, G., G\"otze, F.: { The Arithmetic of Distributions in Free Probability Theory}. {
Preprint}, Arxiv, {http://www.arxiv.org/abs/math.OA/0508245}.

\bibitem{CollingwoodL} Collingwood, E. F., Lohwater, A. J.: {
The theory of cluster sets}. {Cambridge Tracts in Mathematics and Mathematical 
Physics}, No. 56 Cambridge University Press, Cambridge {(1966)} 



\bibitem{Denjoy} Denjoy, A.: {Sur l'iteration des fonctions analytiques}. C. R.
Acad. Sci. Paris {\bf 182},
255--257 (1926).



\bibitem{Garnett} Garnett, J. B.: {Bounded analytic functions}. Academic Press, New
York, (1981).



\bibitem{Nev} Nevanlinna, R.: {Analytic functions}. {Translated from the second German edition by 
Philip Emig. Die Grundlehren der mathematischen Wissenschaften, Band 162},
 {Springer-Verlag}, {New York}, {(1970)}




 \bibitem{Shapiro} Shapiro, J. H.: {Composition operators and classical function
theory}.
 Springer, New York, (1993).

\bibitem{SteinWeiss} {Stein, E. M., Weiss, G.}: {Introduction to {F}ourier analysis 
on {E}uclidean spaces}. {Princeton Mathematical Series, No. 32},
{Princeton University Press},
{Princeton, N.J.},
{(1971)}.



\bibitem{Voiculescu1} {Voiculescu, D.}: {Addition of certain noncommuting random
variables}.
{J. Funct. Anal.} {\bf 66}, {323--346} {(1986)}.


\bibitem{V3} {Voiculescu, D.}: {The analogues of entropy and of Fisher's
information measure in free probability theory. I}. {Comm. Math. Phys.}
 {\bf 155},  {411--440}  {(1993)}.



\bibitem{VDN} {Voiculescu, D. V., Dykema, K. J., Nica, A.}: {Free Random
Variables}. {CRM Monograph Series, Vol. 1}, {American Mathematical Society,
Providence, RI, (1992)}.

\bibitem{Wolff} {Wolff, J.}: {Sur l'iteration des fonctions born\'ees}. C. R.
Acad. Sci. Paris {\bf182},
200-201 (1926).

\end{thebibliography}
\end{document}